\documentclass{article}

\usepackage{amsmath,amsfonts}

\usepackage[OT2,T1]{fontenc}
\DeclareSymbolFont{cyrletters}{OT2}{wncyr}{m}{n}
\DeclareMathSymbol{\Sha}{\mathalpha}{cyrletters}{"58}

\author{Brian Lawrence}
\title{Nonvanishing of Central Derivatives of Modular $L$-series in Level $p^2$}

\begin{document}

\maketitle

        A quadratic twist of the L-function associated with a modular form is
known to satisfy a functional equation, which may be even or odd.  A
result due to Gross and Zagier explicitly computes the central value of
the L-function or its derivative.  In prime level when the functional equation is even,
Michel and Ramakrishnan have used an averaging method to prove several
consequences of the Gross-Zagier formulae, including a non-vanishing
result.  The present research concerns L-functions arising from newforms in prime-squared level, which necessarily have odd functional
equations.  Such an L-function has a central value of zero;
the Gross-Zagier formulae compute the central value of its derivative. 
Using the Michel-Ramakrishnan averaging method, we compute the average value of these derivatives over different
L-functions.  In
particular, we show that under suitable conditions there exists an
L-function with only a simple zero at the center of symmetry.
The proof requires us to bound the contribution of oldforms from level $p$.

\section{Introduction}

Let $f$ be a cusp form of weight $2k$ and level $N$.  Let $K$ be an imaginary quadratic field with odd discriminant $D$ and ring of integers $R$, and let $\chi$ be a character of the ideal class group of $R$.  In \cite{GZ}, Gross and Zagier study the $L$-function of the twisted modular form $f \times \Theta_\chi$.  In particular, they show that this $L$-function satisifies a functional equation with center at $s=k$, with sign determined by the splitting of the prime factors of $N$ in $D$.  Furthermore, Gross and Zagier determine the leading behavior of the $L$-function at $k$ -- the value $L(k, f \times \Theta_\chi)$ if the equation is even and the derivative $L'(k, f \times \Theta_\chi)$ if it is odd -- as a Petersson inner product $(f, g)$, where $g$ is itself a cusp form of weight $2k$ and level $N$.

The aim of this paper is to prove a nonvanishing result for the derivative $L'(k, f \times \Theta_\chi)$.  Specifically, suppose $N = p^2$ with $p$ prime.  Then the functional equation satisfied by the twisted $L$-function is necessarily odd, so the value $L(k, f \times \Theta_\chi)$ is zero.  We wish to show that, for sufficiently large $p$, there is a newform $f$ in level $p^2$ such that the $L$-function has only a single zero at $k$.

Interest in non-vanishing of derivatives of modular $L$-functions dates to work of Kolyvagin  \cite{Koly}, proving finiteness of the Tate--Shafarevich group for modular elliptic curves, assuming some quadratic twist of the $L$-function has a simple zero at $s = 1$. 
Kolyvagin's assumption was proven independently by Murty--Murty \cite{MM} and Bump--Friedberg--Hoffstein \cite{BFH};
Iwaniec \cite{Iwa} refined their nonvanishing result to a more precise quantitative estimate.

Our approach is based on the averageing method of Michel and Ramakrishnan \cite{MR}.  This method has its origins in a calculation of Gross and Zagier \cite{GZ}, based on the observation that any functional on a space of modular forms -- say, the functional that assigns to $f$ a central $L$-value -- can be represented as a Petersson inner product of $f$ against a fixed kernel.  Michel and Ramakrishnan observed that this method allows one to compute the average central value over cusp forms $f$; their method was extended by Goldfeld--Zhang \cite{GolZh}, Nelson \cite{Nel} and others.

First we will extend the calculation in \cite{GZ} to find a holomorphic cusp form $g$  of level $N$ and weight $2k$ such that, for all $f$ of level $p^2$ and weight $2k$, the central value of the derivative is given by the Petersson inner product
$$L'(k, f \times \Theta_\chi) = (f, g).$$ 
Given an orthogonal (but not necessarily orthonormal) basis of functions $f_i$, we can expand $g$ in terms of this basis as
$$\sum \frac{(f_i,g)}{(f_i,f_i)} f_i = g.$$
Writing the inner product $(f_i, g)$ as the derivative of an $L$-function, we get
$$\sum \frac {L'(k, f_i \times \Theta_\chi)} {(f_i,f_i)} f_i = g.$$
Finally, we can apply a linear functional to this equality of forms, say by taking the first Fourier coefficient of each side.  We obtain
\begin{equation}
\label{eqn:michelram}
\sum \frac {L'(k, f_i \times \Theta_\chi)} {(f_i,f_i)} a_1(f_i) = a_1(g).
\end{equation}
This is the ``averaging technique'' of Michel and Ramakrishnan.

In order to finish the problem, we will determine the asymptotic behavior of $a_1(g)$ for large $p$.  We will show that for sufficiently large $p$, the contribution from old forms $f$ to the left-hand side of Equation (\ref{eqn:michelram}) is asymptotically smaller than $a_1(g)$.  Thus, some newform $f$ must contribute a nonzero term.

Our methods will prove an unconditional nonvanishing result only for primes $p$ that remain inert in the quadratic field $K$.  For split primes, we prove a nonvanishing result conditional on the nonnegativity of the central derivative $L' (k, f \times \Theta_\chi)$.  Furthermore, in both cases, the effective bound we obtain for $p$ in the nonvanishing result will be conditional on this conjectured nonnegativity.

\section{Extending the Gross-Zagier Formulas}

To begin with we generalize the calculation in Chapter 4 (starting from p. 267) of \cite{GZ} to give results valid for all forms, old and new.

In the discussion leading to Proposition 4.1.2 on page 271, Gross and Zagier show that for any form $f$ of level $N$ and any ideal class $A$,
$$(4 \pi) ^{-s-2k+1} \Gamma(s+2k-1) L_A (f, s+2k-1) = (f, Tr^{ND}_N \Theta_A E_s),$$
where $E_s(z)$ is given by
\begin{equation}
\label{eqn:1}
E_s(z) = \sum_{e|N} \frac{\mu(e)\epsilon(e)}{e^{2s+2k-1}} (N/e)^{-s} E_s^{(1)} \left( \frac{N}{e} z \right).
\end{equation}
Here $E_s^{(1)}$ is the Eisenstein series
$$\frac{1}{2} \sum_{c, d \in \mathbb{Z} \\ D | c} \frac {\epsilon(d)} {(cz+d) ^{2k-1}} 
\frac{y^s} {\left| cz+d \right| ^{2s}}.$$

In the computation of $Tr^{ND}_N \Theta_A E_s$ that follows, Gross and Zagier discard the terms $e<N$, because they contribute only an old form, which is orthogonal to any new form $f$.  Since we need a result valid for all forms $f$, we repeat their calculation, keeping all terms of the sum for $E_s$.

Thus, we begin with a modified Proposition 4.1.2.  If, for all $\nu|N$, we write
$$\tilde{\Phi}_{s,\nu}(z)= Tr^{ND}_N \Theta_A(z) E_s^{(1)} (\nu z)$$
and
\begin{equation}
\label{eqn:phi}
\tilde{\Phi}_{s,mod}(z) = \sum_{e|N} \frac {\mu(e)\epsilon(e)} {e^{2s+2k-1}} (N/e)^{-s} \tilde{\Phi}_{s,N/e}(z),
\end{equation}
then we have
\begin{equation}
\label{eqn:lofs}
(4 \pi) ^{-s-2k+1} \Gamma (s+2k-1) L_A(f,s+2k-1) = (f, \tilde{\Phi}_{s,mod}(z)).
\end{equation}
The correspondence between Gross-Zagier's notation and ours is: Gross-Zagier's $\tilde{\Phi}_s(z)$ is our $\tilde{\Phi}_{s,N}(z)$, up to a factor of $N^k$.

By a calculation analogous to that in \cite{GZ}, we can compute the Fourier coefficients of $\tilde{\Phi} _{s, \nu}$. Then we differentiate at $s = k$ and apply the holomorphic projection lemma from page 292. We omit the details of the calculation.

When $\epsilon(\nu) = 1$, we find that the holomorphic projection is
\begin{eqnarray*}
\Phi_\nu(z) & = & \sum_{m=1}^{\infty} 2^{2k-1} \pi^k \left| D \right| ^{-\frac{1}{2}} \nu^{-k+1} \frac{(k-1)!} {(2k-2)!} m^{k-1} \\
& & \Big[ -\sum_ {0 < n \leq \frac{m \left| D \right|} {\nu}} -P_{k-1} \left(1 - \frac{2 n \nu}{m \left| D \right|} \right) \sigma_\nu'(n) r_A (-n \nu + m \left| D \right|) \\
& & - \sum_{n=0}^{\infty} \left( -2Q_{k-1} \left(1 + \frac{2 n \nu}{m \left| D \right|} \right) \sigma_\nu'(n) r_A (n \nu + m \left| D \right|) \right) \\
& & + \frac{h}{u} r_A(m) \left(\log \frac{\nu \left| D \right| y}{m} - 2 \log{2 \pi} + 2 \frac {\Gamma'} {\Gamma}(k) + 2 \frac{L'}{L}(1, \epsilon) \right) \Big] e^{2 \pi imz}.
\end{eqnarray*}

In case $\epsilon(\nu) = -1$, the $L$-function has even functional equation. In this case Gross and Zagier compute only its value, not its derivative, at the center of symmetry. But the even functional equation enables us to express the derivative in terms of the value of the function. Specifically, using functional equation $4.4.10$ of \cite{GZ}, we see that the function $e(n,y)$ satisfies
\begin{equation}
\label{eqn:diff}
\frac{e^\ast (n,y)} {e_{1-k}(n,y)} = \log  \pi - \log \left| D \right|  - \frac {\Gamma'} {\Gamma} (k).
\end{equation}
Using this result, the rest of the calculation proceeds as in \cite{GZ}.

After holomorphic projection, we have
$$\Phi_\nu (z) = \left( 2 \log  2 \pi - \log \frac{N^2}{\nu} \left| D \right|  - 2 \frac {\Gamma'} {\Gamma} (k) \right) \Phi _{1-k, \nu} (z),$$
where
\begin{eqnarray*}
\Phi _{1-k, \nu} (z) & = & \sum_{m=1}^{\infty} 2^{2k-1} \pi^k \left| D \right| ^ {-\frac{1}{2}} \nu ^{-k+1} \frac{(k-1)!} {(2k-2)!} m^{k-1} \\
& & \left[ \frac{h}{u} r_A(m) + \sum _{0 < n \leq \frac{m \left| D \right|} {\nu}} P_{k-1} ( 1 - \frac{2 n \nu} {m \left| D \right|}) \sigma _{\nu, A} (n) r_A (m \left| D \right| - n \nu) \right],
\end{eqnarray*}

Thus, for any holomorphic cusp form $f$ of weight $2k$ and level $N$, the derivative of the twisted $L$-function at $k$ is given by the Petersson inner product
$$L'_A(f,k) = (f, g),$$
where
$$g(z) = \frac{ (4 \pi) ^{k}} { (k-1)! } \sum _{e | N} \mu(e) \epsilon(e) N^{k-1} e^{-k} \Phi_{N/e}^\ast (z);$$
and the value is given by
$$L_A (f, k) = \frac{ (4 \pi) ^k} { (k-1)!} \sum _{e | N} \mu(e) \epsilon(e) N^{k-1} e^{-k} \left(f, \Phi_{1-k, N/e}^\ast (z) \right).$$

\section{Asymptotic Estimates}
\label{sec:asymptoticestimates}

Here we prove symptotic estimates for the Fourier coefficients of the form $\Phi_\nu$ arising in our modified Gross-Zagier calculations, for large $\nu$.  (The first section is from last year.)

\subsection{ $\epsilon(\nu) = 1$ , $\nu$ large}
\label{sec:asymptotics}

Write
$$\Phi_\nu(z) = 2^{2k-1} \pi^k \left| D \right| ^{-\frac{1}{2}} \nu^{-k+1} \frac{(k-1)!} {(2k-2)!} m^{k-1} a_{m,\nu} e^{2 \pi imz},$$
where
$$a_{m, \nu} = \sum_{n \leq \frac{m \left| D \right|}{\nu}} a_{m, n, \nu}.$$
Here $a_{m, n, \nu}$ is given by
$$a_{m, n, \nu} = -P_{k-1} \left(1 - \frac{2 n \nu}{m \left| D \right|} \right) \sigma_\nu'(n) r_A (-n \nu + m \left| D \right|)$$
for $n>0$ (and $n < \frac{m \left| D \right|}{\nu}$),
$$a_{m, -n, \nu} = 2Q_{k-1} \left(1 + \frac{2 n \nu}{m \left| D \right|} \right) \sigma_\nu'(n) r_A (n \nu + m \left| D \right|)$$
for $-n<0$, and
$$a_{m, 0, \nu} = \frac{h}{u} r_A(m) \left( \log \frac{\nu \left| D \right|}{m} - 2 \log{2 \pi} + 2 \frac {\Gamma'} {\Gamma}(k) + 2 \frac{L'}{L}(1, \epsilon) \right).$$

For $\nu > m \left| D \right|$ there are no terms $a_{m, n, \nu}$ with $n>0$.

Now we claim that the sum of the negative terms, $\sum_{n=1}^{\infty} a_{m, -n, \nu}$, grows as $o(1)$ for large $\nu$.  If we can show this, it will follow that the coefficient $a_{m, \nu}$ is dominated by the $\log \nu$ term coming from $a_{m, 0, \nu}$.

A straightforward calculation shows that
$$\left| a_{m, -n, \nu} \right| < 2^{8-2k}n^{-k+1/2}\nu^{-k+1/4}(m \left| D \right|)^k.$$
Summing over $n$, we see that the sum is less than
$$2^{10-2k}\nu^{-k+1/4}(m \left| D \right|)^k$$
in absolute value.  In particular, this sum is bounded in $\nu$.

Thus, we have
$$a_{m, \nu} = \frac{h}{u} r_A(m) \left(\log \frac{\nu \left| D \right|}{m} - 2 \log{2 \pi} + 2 \frac {\Gamma'} {\Gamma}(k) + 2 \frac{L'}{L}(1, \epsilon) \right) + o(1),$$
where the error term is bounded by 
$$2^{10-2k}\nu^{-k+1/4}(m \left| D \right|)^k,$$
for $\nu > m \left| D \right|$.

\subsection{ $\epsilon(\nu) = -1$}
Again, write
\begin{eqnarray*}
\Phi_\nu(z) & = & \left( 2 \log  2 \pi - \log \frac {N^2} {\nu} \left| D \right|  - 2 \frac {\Gamma'} {\Gamma} (k) \right) \\
& & \left[ \sum_{m=1}^{\infty} 2^{2k-1} \pi^k \left| D \right| ^{-\frac{1}{2}} \nu^{-k+1} \frac{(k-1)!}{(2k-2)!} m^{k-1} a_{m,\nu} \right] e^{2 \pi imz},
\end{eqnarray*}
where
$$a_{n, \nu} = \sum_{0 \leq n \leq \frac{m \left| D \right|}{\nu}} a_{m, n, \nu}.$$
Here $a_{m,n,\nu}$ is zero unless $0 \leq n \leq \frac{m \left| D \right|}{\nu}$.

Thus, for $\nu > m \left| D \right|$ we have exactly
$$a_{n, \nu} = \frac{h}{u} r_A(m).$$

\subsection { $\nu$ = 1 }

For $0 < n \leq m \left| D \right|$ we have
$$a_{m, n, \nu=1} = -P_{k-1} \left(1 - \frac{2 n}{m \left| D \right|} \right) \sigma_\nu'(n) r_A (-n + m \left| D \right|).$$

Using elementary estimates for $\sigma_A(n)$, $r_A(n)$ and the Legendre polynomials, we find that for $0 < n \leq m \left| D \right|$,
$$\left| a_{m, -n, 1} \right| \leq 81 2^\frac{-3}{4} (m \left| D \right|) ^\frac{1}{2} (\log m \left| D \right| + 2).$$
If $n > m \left| D \right|$ then
$$\left| a_{m, -n, 1} \right| \leq \frac{9}{8} 2^{1-k} 81 n^\frac{1}{4} (n + m \left| D \right|) ^\frac{1}{4} \left(1 + 2 \frac{n} {m \left| D \right|} \right)^{-k}.$$
Summing over all positive $n$, we find after some simplification that
$$\left| \sum_{n > 0} a_{m, -n, 1} \right| \leq 2^6 ( m \left| D \right| ) ^\frac{3}{2} (\log m \left| D \right| + 3).$$

For $n=0$ we have the equality
$$a_{m, 0, 1} = \frac{h}{u} r_A(m) \left(\log \frac{\left| D \right|}{m} - 2 \log{2 \pi} + 2 \frac {\Gamma'} {\Gamma}(k) + 2 \frac{L'}{L}(1, \epsilon) \right).$$
Summing over all $n$, we have
$$\left| a_{m, 1} \right| \leq 192 ( m \left| D \right| ) ^\frac{3}{2} (\log m \left| D \right| + 1) + \frac{h}{u} r_A(m) \left| (\log \frac{\left| D \right|}{m} - 2 \log{2 \pi} + 2 \frac {\Gamma'} {\Gamma}(k) + 2 \frac{L'}{L}(1, \epsilon) ) \right|.$$

\subsection{A Bound for the Dirichlet $L$-Function}
We will use the estimate
$$\left| \frac{L'}{L} (1, \epsilon) \right| \leq \log \left| D \right|.$$

\subsection{Final Asymptotics of Gross-Zagier}

For our purposes it will be convenient to have large-$N$ asymptotic estimates for the Fourier coefficients of $g$.  From our asymptotic work we know that for $\nu > m \left| D \right|$, if $\epsilon (\nu) = 1$, then the $m$-th Fourier coefficient of $\Phi_\nu$ behaves as 
\begin{eqnarray*}
& & 2^{2k-1} \pi^{k} \frac {(k-1)!} {(2k-2)!} \left| D \right| ^{-\frac{1}{2}} \nu^{-k+1} m^{k-1} \\
& & \left[ \frac{h}{u} r_A(m) \left(\log \frac{\nu \left| D \right|}{m} - 2 \log{2 \pi} + 2 \frac {\Gamma'} {\Gamma}(k) + 2 \frac{L'}{L}(1, \epsilon) \right) + o(1) \right],
\end{eqnarray*}
where the error term is bounded by 
$$ 2^{10-2k}\nu^{-k+1/4}(m \left| D \right|)^k.$$

If $\epsilon (\nu) = -1$ then the $m$-th Fourier coefficient of $\Phi_ \nu$ equals exactly
$$\left( 2 \log  2 \pi - \log \frac {N^2} {\nu} \left| D \right|  - 2 \frac {\Gamma'} {\Gamma} (k) \right) 2^{2k-1} \pi^k \left| D \right| ^{-\frac{1}{2}} \nu^{-k+1} \frac{(k-1)!}{(2k-2)!} m^{k-1} \frac{h}{u} r_A(m).$$

Now suppose $N = p^2$, so that from Equation (\ref{eqn:1}) we have
$$g = \frac{(4 \pi) ^k} {(k-1)!} \left( p^{2k-2} \Phi_{p^2} - \epsilon(p) p^{k-2} \Phi_{p} \right).$$
Using the estimates above we find
$$a_m(g) = 2^{4k-1} \pi^{2k} \left| D \right| ^{-\frac{1}{2}} m^{k-1} \left(\frac{h}{u} r_A(m) \log {p^2} + O(1) \right),$$
regardless of the sign of $\epsilon(p)$.

More precisely, if $\epsilon(p) = 1$ then we have 
\begin{eqnarray*}
a_m(g) & = & 2^{4k-1} \pi^{2k} \left| D \right| ^{-\frac{1}{2}} m^{k-1} \Big[ \frac{h}{u} r_A(m) ( (2 - p^{-1}) \log p \\
& & + (1 - p^{-1}) (\log \frac{\left| D \right|}{m} - 2 \log{2 \pi} + 2 \frac {\Gamma'} {\Gamma}(k) + 2 \frac{L'}{L}(1, \epsilon) ) + o(1)) \Big]\\.
\end{eqnarray*}
Here the error term $o(1)$ is bounded by
$$2^{10-2k}(m \left| D \right|)^k (p^{-2k+1/2} + p^{-k-3/4}).$$

If instead $\epsilon(p) = -1$, then
\begin{eqnarray*}
a_m(g) & = & 2^{4k-1} \pi^{2k} \left| D \right| ^{-\frac{1}{2}} m^{k-1} \\
& & \Big[ \frac{h}{u} r_A(m) ( (2 - 3 p^{-1}) \log p + (1 - p^{-1}) \log \left| D \right| - \log m \\
& & - 2 (1 - p^{-1}) \log{2 \pi} + 2 (1 - p^{-1}) \frac {\Gamma'} {\Gamma}(k) + 2 \frac{L'}{L}(1, \epsilon)) + o(1) \Big],
\end{eqnarray*}
with the error term bounded by
$$ 2^{10-2k}p^{-2k+1/2}(m \left| D \right|)^k.$$

If $N=p$ is a prime then
$$g = \frac{(4 \pi) ^k} {(k-1)!} (p^{k-1} \Phi_p - \epsilon(p) p^{-1} \Phi_1),$$
so if $p$ splits in $K$, then
$$a_m(g) = \frac {2^{4k-1} \pi^{2k}} {(2k-2)!} \left| D \right| ^{-\frac{1}{2}} m^{k-1} (\frac{h}{u} r_A(m) \log {p} + O(1)).$$
Specifically,
\begin{eqnarray*}
a_m(g) & = & \frac {2^{4k-1} \pi^{2k}} {(2k-2)!} \left| D \right| ^{-\frac{1}{2}} m^{k-1} \\ & & \Big[ \frac{h}{u} r_A(m) \big(\log p + (1 - p^{-1}) \log \frac{\left| D \right|}{m} - 2 (1 - p^{-1}) \log{2 \pi} \\
& & + 2 (1 - p^{-1}) \frac {\Gamma'} {\Gamma}(k) + 2 (1 - p^{-1}) \frac{L'}{L}(1, \epsilon) \big) + o(1) \Big],
\end{eqnarray*}
and the error term is bounded by
$$ 2^{10-2k}p^{-k+1/4}(m \left| D \right|)^k + 192 ( m \left| D \right| ) ^\frac{3}{2} (\log m \left| D \right| + 1) p^{-1}.$$

If instead $N=p$ is inert in $K$ then
\begin{eqnarray*}
a_m(g) & = & \frac {2^{4k-1} \pi^{2k}} {(2k-2)!} \left| D \right| ^{-\frac{1}{2}} m^{k-1} \\
& & \Big[  \frac {h} {u} r_A(m) \big(2 (1 - p^{-1}) \log 2 \pi - \log p - (1 - p^{-1}) \log \left| D \right| \\
& & - 2 (1 - p^{-1}) \frac {\Gamma'} {\Gamma} (k) + 2 p^{-1} \frac{L'}{L} (1, \epsilon) \big) + o(1) \Big],\\
\end{eqnarray*}
with error term bounded by
$$192 ( m \left| D \right| ) ^\frac{3}{2} (\log m \left| D \right| + 1) p^{-1}.$$

Finally, if $N=1$ then 
$$g = \frac{(4 \pi) ^k} {(k-1)!} \Phi_1$$
so
\begin{eqnarray*}
a_m(g) & = & \frac{ 2^{4k-1} \pi^{2k}} {(2k-2)!} \left| D \right| ^{-\frac{1}{2}} m^{k-1} \frac{h}{u} r_A(m) \\
& &  \left(\log \frac{\left| D \right|}{m} - 2 \log{2 \pi} + 2 \frac {\Gamma'} {\Gamma}(k) + 2 \frac{L'}{L}(1, \epsilon) + E \right),
\end{eqnarray*}
where the error term $E$ is bounded by
$$\left| E \right| \leq 192 ( m \left| D \right| ) ^\frac{3}{2} (\log m \left| D \right| + 1).$$

\section{Miscellaneous Calculations}

To apply the Michel-Ramakrishnan averageing formula we will need an orthogonal basis for the space of old forms, and estimates of their $L$-functions.  Toward this end, we collect here some results on the $L$-functions and inner products of old forms.

\subsection{Computing $L$-functions}
\label{sec:Lfunc}

For old forms in level $p^2$, which come from levels $1$ and $p$, we need to express the $L$-function in terms of the $L$-functions in lower levels.

Suppose first that $g$ is a (new) Hecke eigenform in level $1$, which gives the three old forms
\begin{eqnarray*}
g_1(z) & = & g(z) \\
g_p(z) & = & g(pz) \\
g_{p^2}(z) & = & g(p^2z) 
\end{eqnarray*}
in level $p^2$.  (If $g$ is new in level $p$ then we define the forms $g_1$ and $g_p$ similarly.)  We need to express the $L$-functions
$$L(s, g_i \times \Theta_A),$$
for $i = 1, p, p^2$, in terms of the $p$-invariant $L$-function
$$L(s, g \times \Theta_A).$$
For $i=1$, the functions $g_1$ and $g$ have the same Fourier coefficients, so their $L$-functions coincide.  For $i=p$ and $i=p^2$, we use the method described in \cite{Shimura}.  Here we provide only a sketch of the calculation due to limitations of space.

Each coefficient of the twisted $L$-function is the product of coefficients of a Hecke $L$-function and the $L$-function associated to a theta function; both of these series have Euler products.  The method in $\cite{Shimura}$ shows how to compute the Euler product of the twisted series given these two Euler products.  It is easy to relate the Hecke $L$-series of $g_i$ to that of $g$.  Using $\cite{Shimura}$ we obtain the results below.

Hence, we have
$$L(s, g_p \times \Theta_\chi) = \frac{(\alpha_p + \beta_p) p^{-s} - a_p(g) p^{-2s}} {1 - p^{2k-1}p^{-2s}} L(s, g \times \Theta_\chi)$$
and
$$L(s, g_{p^2} \times \Theta_\chi) = \frac{ c_2 p^{-2s} - c_3 p^{-3s} + c_4 p^{-4s}}{1 - \alpha_p \beta_p \gamma_p \delta_p p^{-2s}} L(s, g \times \Theta_\chi),$$
where
\begin{eqnarray*}
c_2 & = & (\alpha_p^2 + 1 + \beta_p^2)\\
c_3 & = & (\alpha_p + \beta_p)a_p(g)\\
c_4 & = & p^{2k-1} \\
\end{eqnarray*}
Although the calculation is only valid for $s$ with sufficiently large real part, this identity extends to all $s$ by analytic continuation.

\subsubsection{L-functions from Level $p$}

Now suppose $h$ is a Hecke eigenform, new in level $p$, so $h_1(z)=h(z)$ and $h_p(z)=h(pz)$ are the two corresponding old forms in level $p^2$.  Again, the $L$-functions for $h_1$ and $h$ coincide, but we still have to compute the $L$-function for $h_p$.

By the same method as in the previous section,
$$L(s, g_p \times \Theta_\chi) = ((\alpha_p + \beta_p) p^{-s} - a_p(g) p^{-2s}) L(s, g \times \Theta_\chi).$$
Again, by analytic continuation the result is true for all $s$, though the calculation is only valid for $s$ with large real part.

We are interested in the derivatives at the central point $s=k$.  If the functional equation is odd, so that $L$ itself vanishes at $k$, then differentiation yields
$$L'(k, g_p \times \Theta_\chi) = ((\alpha_p + \beta_p) p^{-k} - a_p(g) p^{-2k}) L(k, g \times \Theta_\chi).$$

If instead the functional equation is even, then we know from logarithmic differentiation of the functional equation that 
$$\frac {L'}{L} (k, g \times \Theta_\chi) = 2 \log 2 \pi - \log p \left| D \right| - 2 \frac {\Gamma'} {\Gamma}.$$
Thus we find 
\begin{eqnarray*}
L' (k, g_p \times \Theta_\chi) & = & L (k, g \times \Theta_\chi) \big[ -\log p ((\alpha_p + \beta_p) p^{-s} - 2 a_p(g) p^{-2s}) \\
 & & + (2 \log 2 \pi - \log p \left| D \right| - 2 \frac {\Gamma'} {\Gamma}(k)) ((\alpha_p + \beta_p) p^{-s} - a_p(g) p^{-2s}) \big]. \\
\end{eqnarray*}

\subsection{Inner Products}

We need to construct an orthonormal basis for the space of forms in level $p^2$.  It is known that the Hecke eigenforms form an ONB for the new forms and that the old forms are orthogonal to the new forms.  But it is still necessary to orthogonalize the standard basis of old forms.

\subsubsection{Inner Products of Old Forms from Level $1$}

To this end, suppose first that $g$ and $h$ are eigenforms in level $1$, so $g_1$, $g_p$, $g_{p^2}$, $h_1$, $h_p$ and $h_{p^2}$ are three associated old forms in level $p^2$.  We wish to calculate inner products of the form $(g_i,h_j)$ in terms of the inner product $(g,h)$ in level $1$. Again, we omit detailed calculations. The general plan of attack is as follows. The inner product can be expressed as an integral over a fundamental domain for $\Gamma_0(p^2)$. Applying the transformation law the integral can be expressed as a sum of integrals over a fundamental domain for $\Gamma_0(SL_2(\mathbb{Z}))$. The sum in turn corresponds to a Hecke trace operator, for which we know that $g$ and $h$ are eigenforms. Writing the traces as scalar multiples of $g$ and $h$, we end up with the inner product expressed in terms of the product $(g,h)$ in level $1$.

Assuming $g$ is normalized, so the eigenvalue of the trace operator $T_p$ is exactly $a_p(g)$, we have
\begin{eqnarray*}
(g_1, h_1) &=& p(p+1) (g,h)\\
(g_p,h_1) &=& p^{-2k+2} a_p(g) (g,h)\\
(g_p, h_p) &=& p(p+1) p^{-2k} (g,h)\\
(g_{p^2},h_1) &=& (p^{2-4k} a_p(g)^2 - p^{1-2k}) (g,h)\\
(g_{p^2},h_p) &=& p^{-4k+2} a_p(g) (g,h)\\
(g_{p^2}, h_{p^2}) &=& p^{1-4k} (p+1) (g,h).\\
\end{eqnarray*}

\subsubsection{Inner Products of Old Forms from Level $p$}

If $g$ and $h$ are new forms in level $p$ then the corresponding forms $g_1, g_p, h_1, h_p$ are old in level $p^2$.  We need to compute the inner products $(g_1, h_1)$, $(g_p, h_1)$ and $(g_p, h_p)$.

Since $X(\Gamma_0(p^2))$ covers $X(\Gamma_0(p))$ exactly $p$ times, the usual transformation shows that
$$(g_1, h_1) = p (g,h)$$
and
$$(g_p, h_p) = p^{1-2k} (g,h).$$

The computation of $(g_p, h_1)$ is more involved; the methods in the previous section give
$$(g_p, h_1) = p^{1-2k} a_p(g) (g,h).$$

\section{Old Form Contributions}

Now we return to the Michel-Ramakrishnan average formula.  Fix a level $N$ and weight $2k$, a quadratic field $K$ of discriminant $D$, and an ideal class $A$.  Let $\epsilon$ be the character corresponding to the discriminant $D$, so for $p$ a prime, $\epsilon(p)=1$ if and only if $p$ splits in $K$.  Suppose now that $p$ is such a prime and $N = p^2$.

For any newform $f$ of level $1$, $p$ or $p^2$, the twisted $L$-function satisfies an odd functional equation.  Thus, we have $L(k, f \times \Theta_\chi) = 0$ for all such newforms, and hence for all forms of level $p^2$.  We are interested in the value of the first derivative $L'(k, f \times \Theta_\chi)$.

We have calculated a holomorphic modular form $g$ of level $N$ and weight $2k$ such that, for all $f$ of level $p^2$ and weight $2k$, the central value of the derivative is given by the Petersson inner product
$$L'(k, f \times \Theta_\chi) = (f, g).$$ 
Given an orthogonal (but not necessarily orthonormal) basis of functions $f_i$, we can expand $g$ in terms of this basis as
$$\sum \frac{(f_i,g)}{(f_i,f_i)} f_i = g.$$
Writing the inner product $(f_i, g)$ as the derivative of an $L$-function, we get
$$\sum \frac {L'(k, f_i \times \Theta_A)} {(f_i,f_i)} f_i = g.$$
Finally, we can apply a linear functional to this equality of forms, say by taking the first Fourier coefficient of each side.  We obtain
$$\sum \frac {L'(k, f_i \times \Theta_A)} {(f_i,f_i)} a_1(f_i) = a_1(g).$$
This is the ``averaging technique'' of Michel and Ramakrishnan.

We wish to prove a nonvanishing result, of the type
$$L'(k, f_i \times \Theta_\chi) \neq 0,$$
for some new form $f_i$ of level $p^2$.  We know from our asymptotic bounds that $a_1 (g)$ is nonzero.  But old forms also contribute to the sum.  Our aim is to show that the old form contribution is asymptotically smaller than $a_1 (g)$.

First choose an orthogonal basis as follows: let $f^{(p^2,i)}$, $f^{(p,i)}$ and $f^{(1,i)}$ denote bases of newforms in levels $p^2$, $p$ and $1$, respectively.  For any form $f$ and integer $q$, let $f_q(z) = f(qz)$.  Then the collection $f^{(p^2,i)}$, $f^{(p,i)}_1$, $f^{(p,i)}_p$, $f^{(1,i)}_1$, $f^{(1,i)}_p$, $f^{(1,i)}_{p^2}$ forms a basis for modular forms in level $p^2$.  This is almost, but not quite, an orthonormal basis, in the following sense: two such forms $f^{({j_1},{i_1})}_{q_1}$ and $f^{({j_2},{i_2})}_{q_2}$ are orthogonal unless $j_1=j_2$ and $i_1 = i_2$.  (This follows from our preliminary calculations.)  Thus, in order to orthogonalize this basis, it suffices to apply Gram-Schmidt separately to every pair $f^{(p,i)}_1$, $f^{(p,i)}_p$, and every triple $f^{(1,i)}_1$, $f^{(1,i)}_p$, $f^{(1,i)}_{p^2}$.

Using the inner products computed above, we find that $f^{(p,i)}_1$ and
$$\tilde{f}^{(p,i)}_p = f^{(p,i)}_p -  p^{-2k} a_p(f^{(p,i)}) f^{(p,i)}_1$$
are orthogonal, as are $f^{(1,i)}_1$,
$$\tilde{f}^{(1,i)}_p = f^{(1,i)}_p - p^{-2k} \frac{p}{p+1} a_p(f^{(p,1)}) f^{(1,i)}_1$$
and
$$\tilde{f}^{(1,i)}_{p^2} = f^{(1,i)}_{p^2} -  C_p  f^{(1,i)}_p -  C_1  f^{(1,i)}_1.$$
Here
$$C_p = \left(1 - \frac{1} {(p+1)^2 - p^{-2k+2} a_p(f)^2} \right) p^{-2k} a_p(f^{(1,i)})$$
and
$$C_1 = \left( \frac{-1} {(p+1)^2 - p^{-2k+2} a_p(f)^2} \right) \frac{p^{-4k+1}} {p+1} a_p (f^{(1,i)})^2 - \frac{p^{-2k}} {p+1}.$$
(The coefficients $C_p$ and $C_1$ are dependent on $f$, but we suppress this dependence in our notation.)
We want to estimate the contribution of each of these terms to the Michel-Ramakrishnan average.

We already have the asymptotic result
$$a_1(g) = 2 \pi \left| D \right| ^{-\frac{1}{2}} r_A(1) \left( \frac{h}{u} \log {p^2} + O(1) \right).$$
Supposing $A$ is the class of princpal ideals, so $r_A(1) = 1$, we have
$$a_1(g) = 2 \pi \left| D \right| ^{-\frac{1}{2}} \left( \frac{h}{u} \log {p^2} + O(1) \right).$$

A word about notation: we use $f^{(j,i)}$ to represent a new eigenform in level $j$, and $f^{(j,i)}_q$ the corresponding old eigenforms in level $p^2$.  Thus, when we write an inner product of the form $(f^{(j,i)},f^{(j,i)})$, we mean for the inner product to be computed in level $j$; but when we write it as $(f^{(j,i)}_q,f^{(j,i)})_q$, we mean the inner product computed in level $p^2$.

\subsection{Old Forms from Level $1$}

We begin with the old forms from level $1$.  Since these forms and their $L$-functions are independent of $p$, we expect that their contribution to the sum will be bounded in $p$.  We now confirm this intuition; in fact, we will see that they decay as $O(p^{-2})$.  (The additional $p^2$ comes from the change of fundamental domain, which increases the $(f,f)$ denominator by a factor of $p(p+1)$.)

For each normalized eigenform $f$ of level $1$ (here we omit the superscript $(1,i)$), let $c_1 = (f,f)$ be the inner product computed in level $1$, and let
$$c_2 = \max_\chi \left| L' (k, f \times \Theta_\chi) \right|,$$
where the maximum is taken over all Dirichlet characters $\chi$ modulo $D$.  We want to show that the three corresponding terms of the Michel-Ramakrishnan average are bounded in terms of $c_1$ and $c_2$.

For $\tilde{f}_1$, since $\Theta_A$ is an average over $\chi$ of the forms $\Theta_\chi$, we have
$$\frac {L'(k, f_1 \times \Theta_A)} {(f_1,f_1)} a_1(f_1) = \frac{c_2} {p(p+1) c_1}.$$

For $\tilde{f}_p$, we need to estimate
$$\frac {L' (k, \tilde{f}_p \times \Theta_A)} { (\tilde{f}_p , \tilde{f}_p)} a_1(\tilde{f}_p),$$
where
$$\tilde{f}_p = f_p - p^{-2k} \frac{p}{p+1} a_p(f) f_1.$$
We can bound the numerator by elementary means. To bound the denominator away from zero we use the Ramanujan bound for $a_p(f)$. Loosely speaking, this bound guarantees that $f_p$ is close to orthogonal to $f_1$, which guarantees that $\tilde{f}_p$ is not too close to zero. Using these bounds, we find that the contribution of $f_p$ to the Michel-Ramakrishnan average is bounded by
$$\left| \frac {L' (k, \tilde{f}_p \times \Theta_A)} { (\tilde{f}_p , \tilde{f}_p)} a_1(\tilde{f}_p) \right| \leq 72 \frac{c_2}{p(p+1)c_1}.$$

By similar means, using elementary bounds and the Ramanujan bound to estimate the $\tilde{f}_p$ contribution, we find that it is bounded as
$$\left| \frac {L' (k, \tilde{f}_{p^2} \times \Theta_A)} { (\tilde{f}_{p^2} , \tilde{f}_{p^2})} a_1(\tilde{f}_{p^2}) \right| \leq \frac{ 32 p^{-2k} c_2} {\frac{1}{10} p^{-4k+2}} (5 p^{-2k-1} c_1) = 1600 p^{-3} \frac{c_2} {c_1}.$$
In short, this bound is $O(p^{-3})$.

Summing over all $f^{(1,i)}_q$, we find that the contribution to the sum from all forms of level $1$ is bounded as $O(p^{-2})$.  In particular, since the dimension of the space of forms of level $1$ and weight $2k$ is at most $\frac{k} {12} + 1$, the total contribution is bounded by 
$$\left( \frac{k} {12} + 1 \right) \frac{c_2}{c_1} ( 73 p^{-2} + 1600 p^{-3}).$$

If we assume further that all the $c_2$ are positive, then by the averageing formula in level $1$ we know that
$$\sum \frac{c_2}{c_1} \leq \frac{ 2^{4k-1} \pi^{2k}} {(2k-2)!} \left| D \right| ^{-\frac{1}{2}} \frac{h}{u} (\log \left| D \right| - 2 \log{2 \pi} + 2 \frac {\Gamma'} {\Gamma}(k) + 2 \frac{L'}{L}(1, \epsilon) + E),$$
where the error term $E$ is bounded by
$$\left| E \right| \leq 192 \left| D \right| ^\frac{3}{2} (\log \left| D \right| + 1).$$
Thus, the contribution to the sum from old forms from level $1$ is bounded by
$$C_k C_p \frac{ 2^{4k-1} \pi^{2k}} {(2k-2)!} \left| D \right| ^{-\frac{1}{2}} \frac{h}{u} \left( \log \left| D \right| - 2 \log{2 \pi} + 2 \frac {\Gamma'} {\Gamma}(k) + 2 \frac{L'}{L}(1, \epsilon) + E \right),$$
where
$$C_k = \frac{k} {12} + 1$$
and
$$C_p = 73 p^{-2} + 1600 p^{-3}.$$

\subsection{Old Forms from Level $p$, where $\epsilon(p) = 1$}

Now we consider the old forms from level $p$.  Their contribution will be bounded by $p^{-1} \log p$.

Let $f$ be a normalized eigenform of level $p$ (here we again suppress the superscript $(1,i)$), let $c_1 = (f,f)$ be the inner product computed in level $p$, and let
$$c_2 = \max_\chi \left| L' (k, f \times \Theta_\chi) \right|,$$
where the maximum is taken over all Dirichlet characters $\chi$ modulo $D$.  We want to estimate the two corresponding terms of the Michel-Ramakrishnan average in terms of $c_1$ and $c_2$. 

First consider the contribution from $f_1$, for each $f$.  We have
$$\frac {L'(k, f_1 \times \Theta_A)} {(f_1,f_1)} a_1(f_1) = \frac{1}{p} \frac {L'(k, f \times \Theta_A)} {(f,f)} a_1(f).$$
By Gross-Zagier in level $p$, the sum of these terms behaves asymptotically as
$$O \left( \frac{\log p}{p} \right).$$
Specifically, the sum is bounded by
\begin{eqnarray*}
p^{-1} \frac {2^{4k-1} \pi^{2k}} {(2k-2)!} \left| D \right| ^{-\frac{1}{2}} \\
\Big[ \frac{h}{u} \big(\log p + (1 - p^{-1}) \log \left| D \right| - 2 (1 - p^{-1}) \log{2 \pi} + 2 (1 - p^{-1}) \frac {\Gamma'} {\Gamma}(k) + 2 (1 - p^{-1}) \frac{L'}{L}(1, \epsilon)  \big) + o(1) \Big],
\end{eqnarray*}
with error term bounded by
$$ 2^{10-2k}p^{-k+1/4}\left| D \right|^k + 192 \left| D \right| ^\frac{3}{2} (\log \left| D \right| + 1) p^{-1}.$$

Next consider the contribution from $\tilde{f}_p$.  To bound this sum we need a positivity result, which we know is unconditionally true: $L'(k, f \times \Theta_A) \geq 0$. (\cite{JN})  Since $a_1(f) = 1$, we have
$$\frac {L'(k, f \times \Theta_A)} {(f,f)} a_1(f) \geq 0$$
for each $f$.  Thus, the Michel-Ramakrishnan bound on the average of these terms translates to individual bounds on each term. Using this result, the usual methods give
$$\left| \frac {L'(k, f_p \times \Theta_A)} {(f_p,f_p)} a_1(f_p) \right| \leq 2 p^{-\frac{3}{2}} \frac{L'(k, f \times \Theta_A)} {(f,f)} a_1(f).$$

Summing over all new forms $f$ of level $p$, and using the asymptotic results to bound the Michel-Ramakrishnan average on the right-hand side, we see that the contribution of all the $f_p$ is at most
\begin{eqnarray*}
& & 2 p^{-\frac{3}{2}} \frac {2^{4k-1} \pi^{2k}} {(2k-2)!} \left| D \right| ^{-\frac{1}{2}} \Big[ \frac{h}{u} \big(\log p + (1 - p^{-1}) \log \left| D \right| \\
& & - 2 (1 - p^{-1}) \log{2 \pi} + 2 (1 - p^{-1}) \frac {\Gamma'} {\Gamma}(k) + 2 (1 - p^{-1}) \frac{L'}{L}(1, \epsilon) \big) + o(1) \Big],
\end{eqnarray*}
where again the error term is bounded by
$$ 2^{10-2k}p^{-k+1/4}\left| D \right|^k + 192 \left| D \right|^\frac{3}{2} (\log \left| D \right| + 1) p^{-1}.$$

\subsection{Old Forms from Level $p$, where $\epsilon(p) = -1$}
The notation is as in the previous section.  As in the case $\epsilon(p) = 1$, we find that
$$\frac {L'(k, f_1 \times \Theta_\chi)} {(f_1,f_1)} a_1(f_1) = \frac{1}{p} \frac {L'(k, f \times \Theta_\chi)} {(f,f)} a_1(f),$$
and this sum behaves asymptotically as
$$O(\frac{\log p}{p}).$$
Precisely, we can estimate the sum as 
\begin{eqnarray*}
& & \frac {2^{4k-1} \pi^{2k}} {(2k-2)!} \left| D \right| ^{-\frac{1}{2}} \Big[ \frac {h} {u} (2 (1 - p^{-1}) \log 2 \pi - \log p \\
& & - (1 - p^{-1})\log \left| D \right| - 2 (1 - p^{-1}) \frac {\Gamma'} {\Gamma} (k) + 2 p^{-1} \frac{L'}{L} (1, \epsilon)) + o(1) \Big],
\end{eqnarray*}
with error term bounded by
$$192 \left| D \right|^\frac{3}{2} (\log \left| D \right| + 1) p^{-1}.$$

By the same methods, we can bound the contribution of $\tilde{f}_p$ by
$$\left| \frac {L'(k, \tilde{f}_p \times \Theta_\chi)} { (\tilde{f}_p, \tilde{f}_p) } a_1 ( \tilde{f} _p) \right| \leq 40 p^{-\frac{3}{2}} (2 \log p + 2 \log 2 \pi k + \log \left| D \right|) \frac {L ( k, f \times \Theta_\chi)} {(f,f)}.$$

Again, the Gross-Zagier estimates show that this is bounded as $o(1)$.  Specifically, the sum of these terms is 
\begin{eqnarray*}
& & 40 p^{-\frac{3}{2}} (2 \log p + 2 \log 2 \pi k + \log \left| D \right|) \frac {2^{4k-1} \pi^{2k}} {(2k-2)!} \left| D \right| ^{-\frac{1}{2}} \\
& & \Big[  \frac {h} {u} \left(2 (1 - p^{-1}) \log 2 \pi - \log p - (1 - p^{-1})\log \left| D \right| - 2 (1 - p^{-1}) \frac {\Gamma'} {\Gamma} (k) + 2 p^{-1} \frac{L'}{L} (1, \epsilon) \right) + o(1) \Big],
\end{eqnarray*}
with error term bounded by
$$192 (\left| D \right| ) ^\frac{3}{2} (\log \left| D \right| + 1) p^{-1}.$$

\section{The Main Result and Conditional Effective Bounds}
We have seen that old forms contribute a term of growth $o(1)$ to the Michel-Ramakrishnan average.  This is true for all primes $p$ if the nonnegativity result holds for the derivative $L'$, and unconditionally true for $p$ with $\epsilon(p) = -1$.  Since the average grows asymptotically as $O(\log p)$, we know that for sufficiently large $p$,
$$\sum_f \frac{L'(k, f \times \Theta_A)} {(f,f)} a_1(f) \neq 0,$$
where the sum is taken over newforms $f$ in level $p^2$ and weight $2k$.  Therefore, there is at least one newform $f$ with
$$L'(k, f \times \Theta_A) \neq 0,$$
conditionally if $\epsilon(p) = 1$ and unconditionally if $\epsilon(p) = -1$.

Using the estimates to the oldform contributions to the Gross-Zagier sum, we can prove the following effective bound, conditionally on the nonnegativity of the central value of $L'$: For all primes $p > 10^4 k \left| D \right|$, there is a newform in level $p^2$ such that the corresponding twisted $L$-function has only a single zero at the central point $k$ of the functional equation.  (This bound is by no means the best possible.)

Using the bounds
$$\log x < 4 x ^ \frac{1}{4} (x>1),$$
$$0 < \frac{\Gamma'}{\Gamma} (k) < \log k,$$
and
$$\left| \frac{L'}{L} (1, \epsilon) \right| < \log \left| D \right|,$$
one easily verifies that the contribution to the Gross-Zagier sum from all old forms combined is less than
$$\frac {2^{4k-1} \pi^{2k}} {(2k-2)!} \left| D \right| ^{-\frac{1}{2}} \frac{h}{u}.$$
But from the result at the end of Section \ref{sec:asymptoticestimates} it is easy to see that the entire Michel-Ramakrishnan sum is necessarily larger than this bound.  Thus, there is a newform that contributes a nonzero value to the sum.

\section{Acknowledgements}
The author would like to thank Dinakar Ramakrishnan for helpful guidance and valuable support.  This work was done at the California Institute of Technology, partially with the support of a Summer Undergraduate Research Fellowship in the summer of 2011, and partially as the author's undergraduate thesis.

\end{document}